\documentclass[14pt]{article}
\usepackage[english]{babel}
\usepackage{amsthm,amsfonts, amsbsy, amssymb,amsmath,graphicx}
\usepackage{graphics}

\newtheorem{thm}{Theorem}

 \newcounter{il}[section]
 \newcounter{il2}[section]

 \newenvironment{dfn}{\trivlist \item[\hskip\labelsep{\bf Definition}]
 \refstepcounter{il}{\bf \arabic{section}.\arabic{il}.}}%
 {\endtrivlist}

 \newenvironment{rk}{\trivlist \item[\hskip\labelsep{\bf Remark}]
 \refstepcounter{il2}{\bf \arabic{section}.\arabic{il2}.}}%
 {\endtrivlist}

 \def\R{{\mathbb R}}
 \def\Z{{\mathbb Z}}
  \def\Q{{\mathbb Q}}

\newcommand{\skcrro}{\raisebox{-0.25\height}{\includegraphics[width=0.5cm]{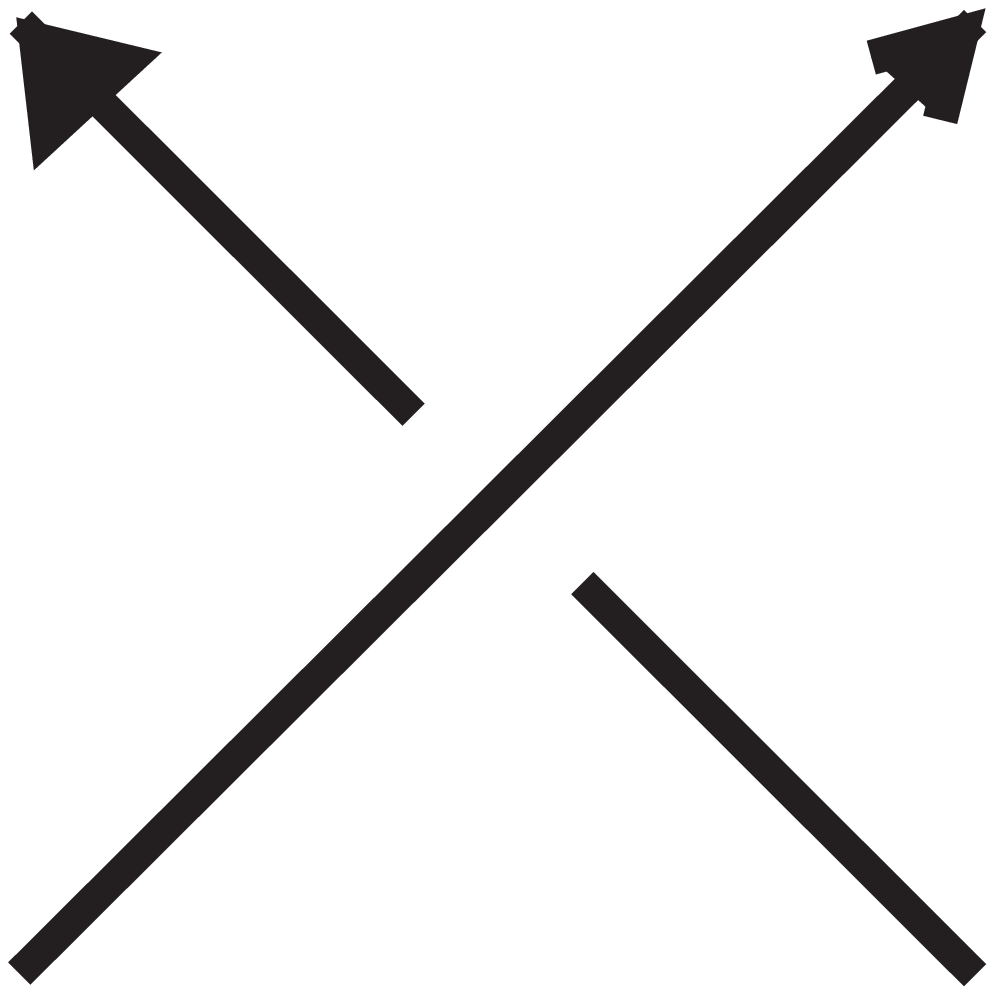}}}
\newcommand{\skcrlo}{\raisebox{-0.25\height}{\includegraphics[width=0.5cm]{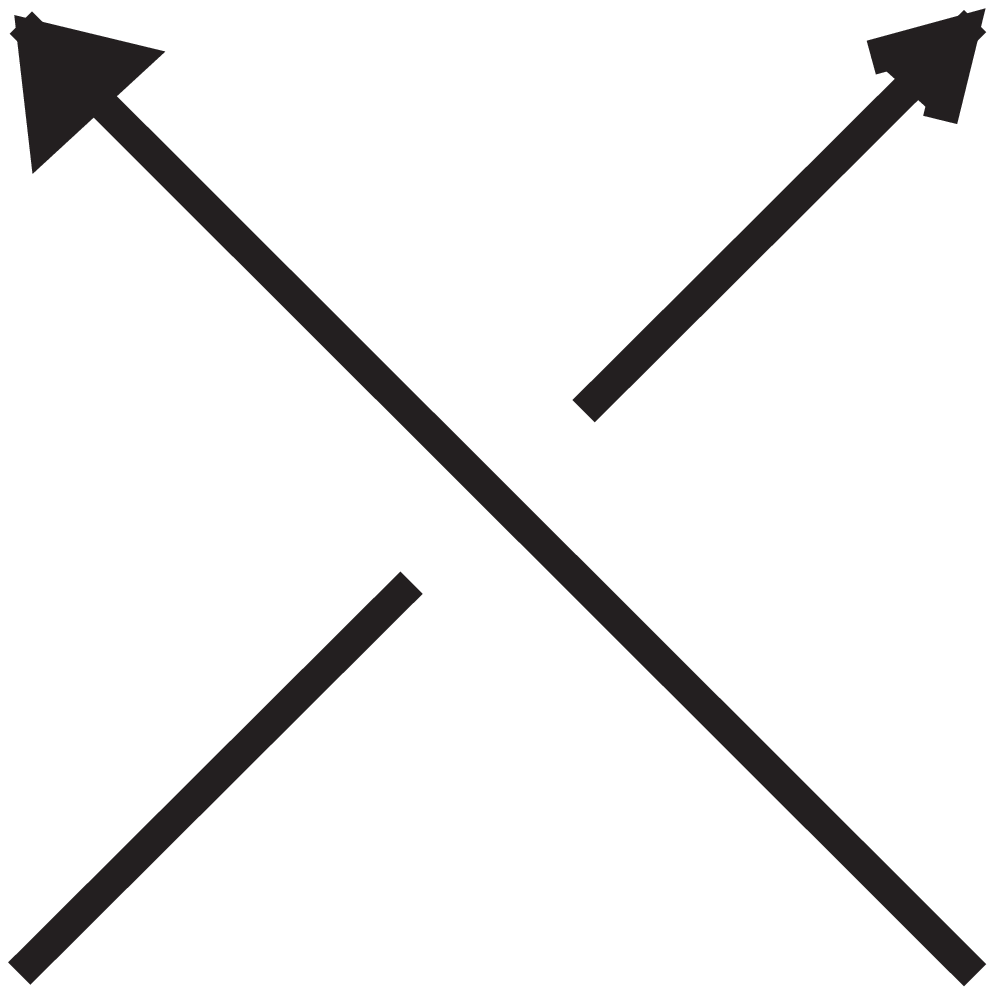}}}
\newcommand{\skcross}{\raisebox{-0.25\height}{\includegraphics[width=0.5cm]{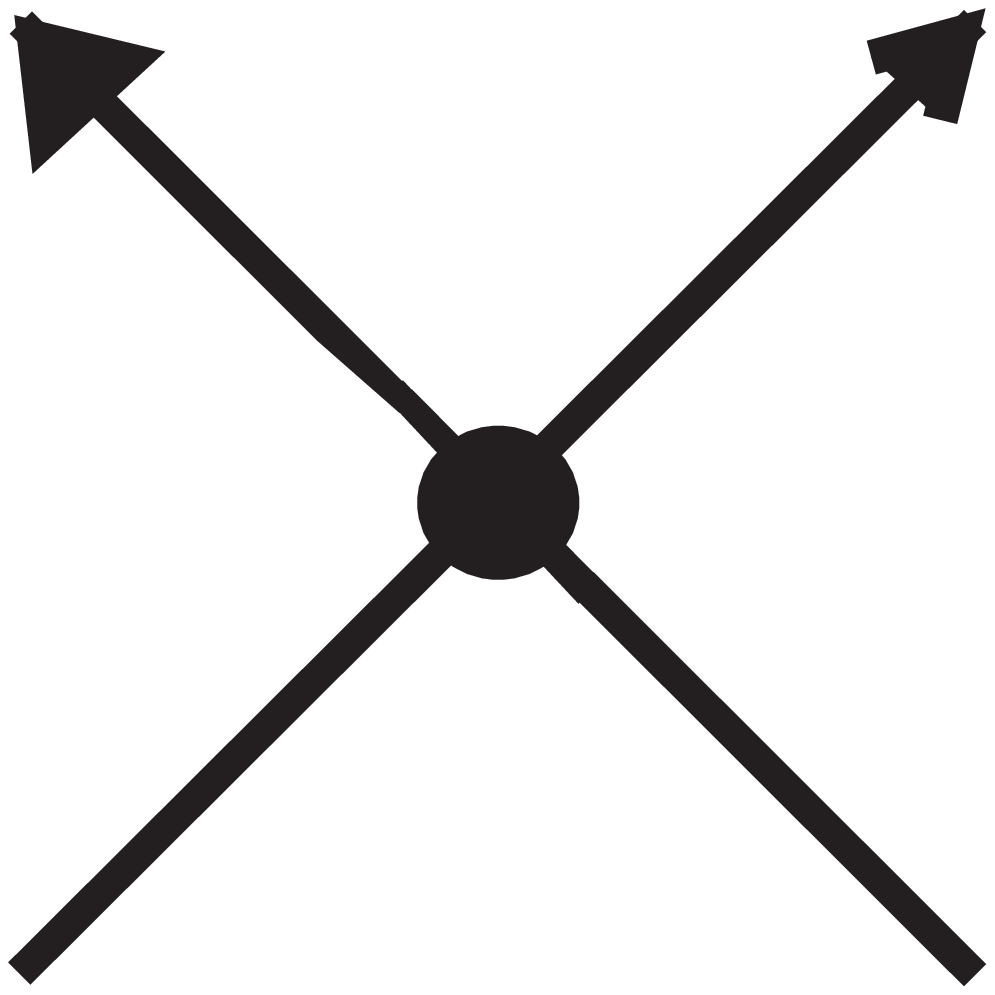}}}

\author{V.O.Manturov\footnote{vomanturov at yandex.ru}}

\title{Free Knots, Groups, and Finite-Type Invariants}

\date{}

\begin{document}

\maketitle

\begin{abstract}
Based on a recently introduced by the author notion of {\em parity},
in the present paper we construct a sequence of invariants (indexed
by natural numbers $m$) of long virtual knots, valued in certain
simply-defined group ${\tilde G}_{m}$ (the Cayley graphs of these
groups are represented by grids in the $(m+1)$-space); the conjugacy
classes of elements of $G_{m}$ play the role of invariants of {\em
compact} virtual knots. By construction, all invariants do not
change under {\em virtualization}. Factoring the group algebra of
the corresponding group by certain polynomial relations leads to
finite order invariants of (long) knots which do not change under
{\em virtualization}.
\end{abstract}

\section{Introduction}

Virtual knot theory was invented by Kauffman, \cite{KaV}; virtual
knots correspond to knots in  thickened surfaces $S_{g}\times
S^{1}$, considered up to isotopies and stabilizations.

The theory of finite type invariants of classical knots was invented
independently by V.A.Vassiliev \cite{Vas} and M.N.Goussarov
\cite{Gus}, and it turned out that many well-known invariants are
expressible in terms of finite-type invariants \cite{BL,BN}; a
breakthrough in Vassiliev's theory of finite type invariants was
marked by the celebrated paper by Kontsevich \cite{Kon}, where the
structure of the space of finite-type Vassiliev invariant was
explicitely described. The Vassiliev knot invariants initiated the
study of virtual knots by Goussarov, Polyak, and Viro \cite{GPV}. In
the latter paper, the authors gave a definition of virtual knots
equivalent to the original one due to Kauffman; this allowed one to
classify invariants of finite type of {\em classical and virtual
knots} in terms of combinatorial formulae, though their definition
of finite type invariants of virtual knots was quite limited. In the
present paper, we shall use a more general (and more natural)
definition of finite-type (Vassiliev) invariants of virtual knots
due to Kauffman \cite{KaV}.

In the case of knots in thickened surfaces (without stabilization),
an analogue of the  Kontsevich theorem was obtained in \cite{AM};
however this solution of the classification problem (an explicit
universal formula) does not look as elegant as in the classical
case; in any case, this formula is quite far away from practical
implementation.

Virtual knots are much more complicated objects than classical ones:
all invariants of classical knots of order zero are constants,
whence the space of finite type invariants of order zero is
infinite-dimensional. The latter statement can be reformulated in
the following manner: there are infinitely many types of virtual
knots where two virtual knots represent the same type whenever one
of them can be obtained from the other by a sequence of
(equivalences and) crossing changes. These ``types'' (or
``equivalence classes'') are called {\em flat virtual knots} (see
\cite{HK}); representing themselves the dual space of the space of
invariants of order $0$ they are important for investigation of
invariants of finite type.

A thorough simplification for the notion of flat virtual knot is the
notion of {\em free knot}. In \cite{Tu}, V.G.Turaev (who first
introduced free knots under the name of ``homotopy classes of Gauss
words'') conjectured that all free virtual knots were trivial; this
conjecture was first disproved by the author in \cite{Ma}, and
later, in \cite{Gib}.

The aim of the present work is to construct invariants of free knots
valued in certain groups $G_{m}$ (depending on a certain natural
parameter $m$), see \cite{MM,MM2}, and to extend this construction
to some simply defined infinite series of finite type invariants for
(long) virtual knots. The groups $G_{m}$ are defined by generators
and relators, and have very simple Cayley graphs; the group $G_{1}$
is isomorphic to the infinite dihedral group.

I am very grateful to O.V.Manturov for fruitful consultations.

\section{Basic Definitions}

Throughout the paper, by {\em graph} we mean a finite multi-graph
(loops and multiple edges are allowed).  Here a $4$-valent graph is
called {\em framed} if for each vertex of it, the four emanating
half-edges are split into two pairs; half-edges belonging to the
same pair are called (formally) {\em opposite}. Assume a framed
$4$-valent graph $K$ is given by its Gaussian diagram $C(K)$.  We
shall use the generic term ``$4$-graph'' to denote a topological
object obtained from a four-valent graph by adding free circles
(without vertices) as connected components.

For framed $4$-graphs one naturally defines unicursal components: if
a graph has connected components homeomeorphic to the circle, they
are treated as unicursal components; for the rest of the graph,
unicursal components are defined as equivalence classes of edges
generated by ``local opposite'' relation at vertices. This relation
of unicursal components naturally yields the {\em number of uncursal
components} which agrees with the number of components of a link
diagram drawn on the plane. Certainly, a framed $4$-graph can be
encoded by a chord diagram if and only if it can be has one
unicursal component.

\begin{dfn}
By a {\em chord diagram} we mean a pair $(S_{0}\sqcup \cdots \sqcup
S_{0}\subset S^{1})$ consisting of an oriented circle (called {\em
the core circle)} and a set of $n$ unordered pairs of points (all
$2n$ points are pairwise distinct). These pairs are called {\em
chords}; points of the pairs are called {\em ends of chords}.
\end{dfn}

 In
figures, we shall depict chords by solid lines connecting pairs of
points on the core circle. We say that a chord diagram is {\em
labeled} if it is marked by a point on the core circle distinct from
chord ends. We consider chord diagrams up to a natural equivalence
(orientation-preseriving homeomorphism of core circle which respects
the collection of chords).

Framed $4$-graphs with one unicursal component are in one-one
correspondence with chord diagrams. The $4$-graph without vertices
corresponds to the chord diagram without chords. Every framed
$4$-valent graph can be represented as the image of the circle going
transversely through all edges. The two preimages of every vertex
correspond to a chord.

One can naturally consider {\em oriented} or {\em non-oriented}
chord diagrams (framed $4$-graphs).

All statements about chord diagrams can be translated into the
language of framed $4$-graphs and vice versa.

By a {\em free knot} \cite{Ma} is meant an equivalence class of
chord diagrams by the equivalence relation generated by the
following three elementary equivalences. The first Reidemeister move
corresponds to an addition/removal of a solitary chord. The second
Reidemeister move is an addition/removal of a couple of {\em
similar} chords. Here chords $a$ and $b$ are called similar if their
chord ends $a_{1},a_{2}$ and $b_{1},b_{2}$ can be numbered in a way
such that $a_{1}$ and $b_{1}$ are adjacent and $a_{2}$ and $b_{2}$
are adjacent (here chord ends $x,y$ are {\em adjacent} if one of the
two components of the complement $C\backslash\{x,y\}$ has no chord
ends inside it).

\begin{figure}
\centering\includegraphics[width=200pt]{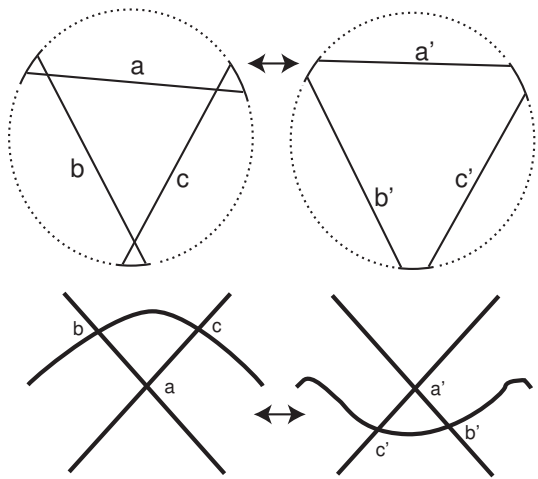} \caption{The
third Reidemeister move and the correspondence between crossings}
\label{sootvtt}
\end{figure}

The third Reidemeister move is depicted in Fig. \ref{sootvtt}.

Here the only changing part of the chord diagram consists of three
segments with two chord ends on each. The chords belonging to the
``stable'' part of the chord diagram are not drawn not depicted, and
the part of the circle containing only stable chord ends is depicted
by dotted lines; in the remaining part, every chord ``moves'' each
of its two ends from one position to the other.

All the three Reidemeister moves on chord diagram originate from the
Reidemeister moves for virtual knots \cite{KaV}, defined later in
the present paper.

Every Reidemeister move transforms one fragment of the frame
four-graph. When depicting such a fragment in a figure we show only
the ``changing part'', leaving the rest of the diagram outside. In
the case of one unicursal component this corresponds to a rule for
transforming a Gauss diagram; this transformation deals with some
arcs of the Gauss diagrams. When depicting on the plane, we do not
show ``stable chords''; we use dotted line for depicting ``stable
parts'' of the core circle.

\begin{figure}
\centering\includegraphics[width=200pt]{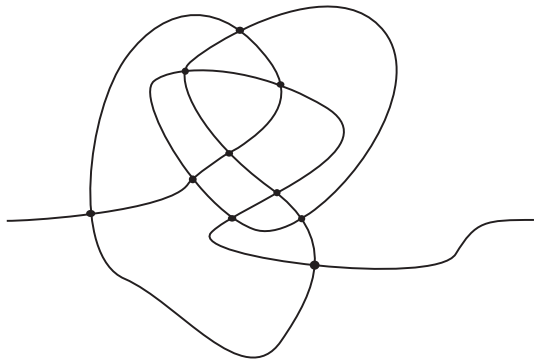} \caption{A
long free knot} \label{lng}
\end{figure}

\begin{dfn}
By a {\em long free knot} we mean an equivalence class of based
chord diagrams by the same Reidemeister move; here we only require
that the marked point lies outside the transformation domain
(equivalently, this marked point lies on a dotted part of the core
circle), see Fig.\ref{lng}.
\end{dfn}

\section{The Groups $G_{m}$. Even and Odd Chords}

Fix a positive integer $m$. Consider the group $G_{m}=\langle
a'_{0},a''_{0},\dots,
a'_{m-1},a''_{m-1},a_{m}|(a_{m})^{2}=e,(a'_{i})^{2}=e,(a''_{i})^{2}=e,i=0,\dots,
m\!-\!1,
a'_{i}a'_{j}=a'_{j}a''_{i},i<j;a'_{i}a''_{j}=a''_{j}a''_{i},i<j;
a'_{i}a_{m}=a_{m}a''_{i} \rangle$.

With every framed chord diagram $D$ we associate an element $G_{m}$.
For this sake, we associate with each chord its {\em index and type}
in the following way.

Let  $f$ be the operation on chord diagrams deleting all odd chords.
It is proved in \cite{Ma} that this operation is well defined: if
two chord diagrams $D,D'$ represent equivalent free knots, then so
are $f(D),f(D')$.

We say that a chord diagram  $D$ a chord has index $0$ if this chord
{\em odd}; those {\em even chords} of $D$ which become odd in
$f(D)$, are decreed to have index $1$, inductively, for $m-1$ we
define those chords of $D$ which remain in $f^{m-1}(D)$ and become
odd ones in $f^{m-1}(D)$, to have index $m-1$; all the remaining
chords are said to have index $m$.

Now, for chords of index $m$, we do not define the type; for chords
of index $k<m$ we define the type ({\em the first type} or {\em the
second type}) depending on the fact whether the number of chords of
the diagram $f^{k}(D)$, linked with the given chord, is {\em even}
(in $f^{k}(D)$) or odd.

\begin{figure}
\centering\includegraphics[width=200pt]{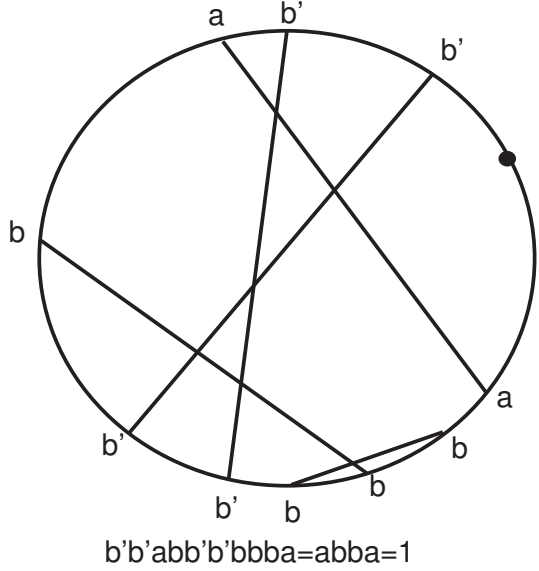}
\caption{Constructing $\gamma(D)$ from a chord diagram}
\label{primer}
\end{figure}

Let us take a chord diagram $D$ and let us start walking from the
base point of the diagram $D$ along the orientation of the circle.
Whenever we meet a chord end of index $i$ and type $a$, we write
down the generator of the group
 $G$ with index $i$ and, if index is less than $m$, with the number of primes equal to
 $a$. If the index is equal to $m$ then we write down the generator
 $a_{m}$.
When we make a full turn returning to the base point, we get a
certain word $\gamma(D)$. Denote the corresponding element of the
group $G$ by $[\gamma(D)]$, see Fig. \ref{primer}.

 In \cite{MM,MM2},
the following statements are proved
\begin{thm}
1. $[\gamma(D)]$ is an invariant of long free knots. 2. The
conjugacy class of the element $[\gamma(D)]$ is an invariant of free
knots.
\end{thm}

The proof of this theorem follows from a straightforward check.

\subsection{The Cayley graph of $G_{m}$}

Elements of $G_{m}$ are in one-one correspondence with integer
points in the Euclidean space ${\bf R}^{m+1}$ with last coordinate
equal to zero or $1$. The unit of the group is represented by the
origin of coordinates. With an element of the group, we associate a
point in ${\bf R}^{m+1}$ defined by induction, as follows. Assume
for some elements $g\in G_{m}$ represented by a word in generators,
we have already defined the corresponding point in the Euclidean
space; let us define those points corresponding to $g\alpha$, where
$\alpha$ is the generator of $G_{m}$. The right multiplication by
the generator with lower index $k$ ($a'_{k}$ or $a''_{k}$ or $a_{k}$
for $k=m$) correspond to the shift of the $(k+1)$-th coordinate. The
direction of the shift is defined as follows: for $0\le k\le m-1$ if
the sum of coordinates
 of the initial point (all except the first $k$) is even then the
 multiplication by $a'_{k}$ increases the first coordinate, whence the multiplication
by  $a''_{k}$ decreases it; if this sum of all coordinates except
the first $k$ is odd, then $a'_{k}$ and $a''_{k}$ change their
roles; the multiplication by $a_{m}$ changes the last  $(m+1)$-st
coordinate: from zero to one and from one to zero. It can be readily
checked that this correspondence is well defined, i.e., the
resulting element of the group $G_{m}$ does not depend on the way of
representing an element of the group as a word in generators.

It can be easily checked that the last coordinated of the element of
$G_{m}$ corresponding to a word coming from a chord diagram, is
zero. The cojugacy class of the element of $G_{m}$ having
coordinates $(x_{0},\dots, x_{m-1})\in \Z^{m}$ (with the last
coordinate equal to zero) consists of elements with coordinates
$(\pm x_{0},\dots, \pm x_{m-1})$ having the last coordinate equal to
zero.

\begin{figure}
\centering\includegraphics[width=350pt]{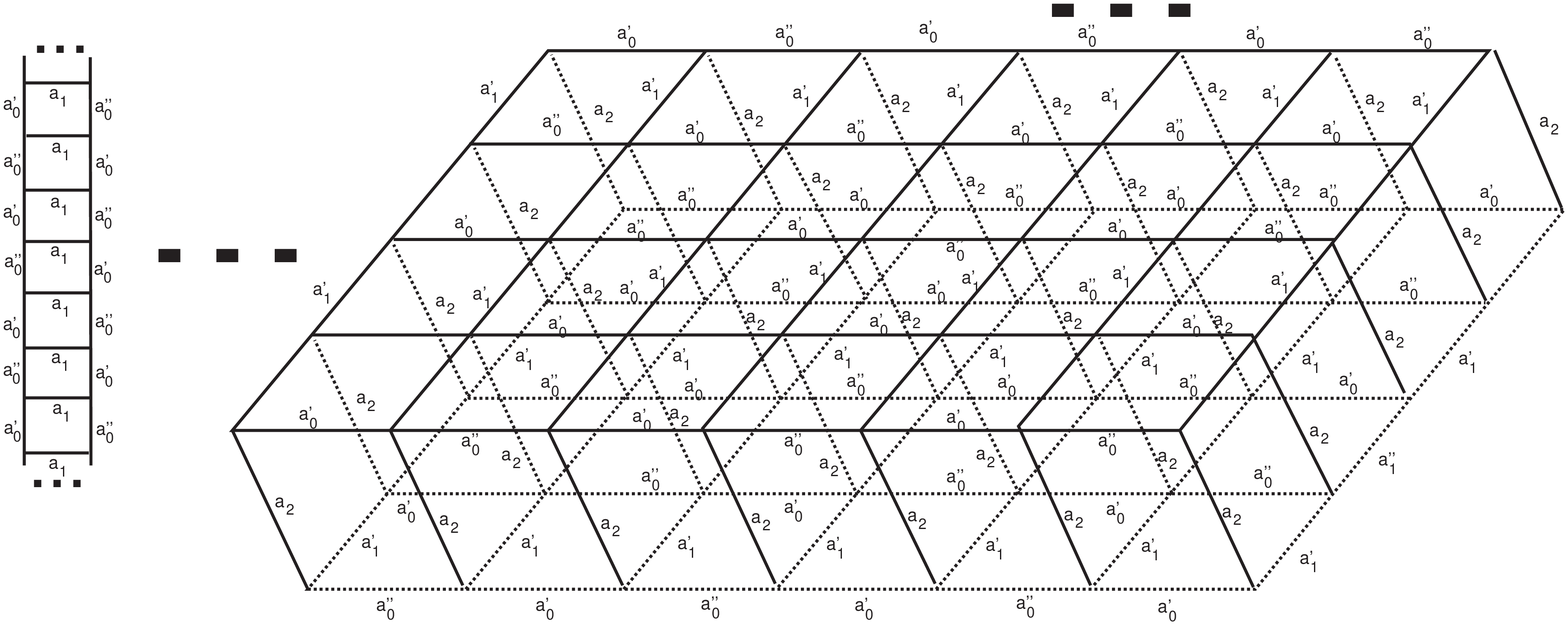} \caption{The
Cayley graphs for $G_1$ and $G_2$} \label{g2}
\end{figure}

The Cayley graphs of $G_{1}$ and $G_{2}$ are given in Fig. \ref{g2}.

\subsection{Virtual knots}

In the present paper, virtual knots appear as a generalization of
free knots with special decorations at vertices of the framed
$4$-graph (resp., chords of the chord diagram). We shall give a
definition in terms of planar diagrams and Reidemeister moves, the
corresponding moves on chord diagram with decorations (called {\em
Gauss diagrams}) can be written straightforwardly.

\begin{figure}
\centering\includegraphics[width=200pt]{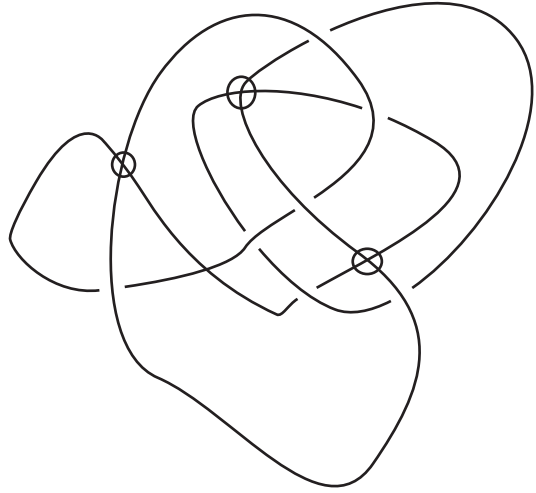} \caption{A
virtual diagram} \label{virt}
\end{figure}

A {\em virtual diagram}, an example see in Fig. \ref{virt}, is a
generic planar immersion of a four-valent graph, where every image
of a vertex is marked as a classical crossing (one pair of opposite
edges is said to form an overcrossing, and the other pair is said to
form an undercrossing); when depicting on the plane, undercrossings
are marked by a broken line; we also encircle the intersection
points between images of edges and say that they from {\em virtual
crossings}. We also admit the case when some components of the
virtual diagram are not graphs, but circles without vertices. For
this sake, we slightly change the notation and will say that a
virtual diagram is a generic immersion of a $4$-graph in $\R^{2}$,
where a $4$-graph means a disjoint sum of a $4$-valent graph with a
collection of circles.

\begin{rk}
Unlike a framed $4$-graph, in a virtual diagram, we have got two new
structures at classical crossings: some pair of opposite edges are
said to form an {\em overcrossing}, the remaining pair forms an {\em
undercrossing}; besides, the four edges incident to classical
crossing obtain a {\em counterclockwise ordering}: being immersed in
the plane, we know not only which edges are opposite, but also which
non-opposite edge follows after the given one in the
counterclockwise direction.
\end{rk}

\begin{dfn}
A {\em virtual link} is an equivalence class of virtual diagram
modulo three classical Reidemeister moves and the detour move, see
\cite{Mybook}. The Reidemeister moves are applied to fragments of
the diagram containing only classical crossings. The detour move
deals with an arc containing only virtual crossings and virtual
self-crossings. This move replaces such an arc with another
(possibliy, self-intersecting) arc having the same endpoints; all
new crossings and self-crossings are marked as virtual ones.
\end{dfn}

For a virtual diagram, one naturally defines {\em unicursal
components}. The immersion image of the $4$-graph is again a
$4$-graph; the latter consists of a $4$-valent graph and a
collection of circles disjoint from the rest of the graph. Every
circle is treated as a unicursal component; other unicursal
components are equivalence classes of edges of $4$-graph components,
where the equivalence is generated by opposite edge relation. This
agrees with the term ``component'' for a link diagram drawn on the
plane. Certainly, the number of unicursal components remains
unchanged under Reidemeister moves.

\begin{dfn}
A {\em virtual knot} is a one-component virtual link. Analogously,
one defines {\em long virtual knots}. In this case, to define {\em
long virtual diagrams}, instead of framed $4$-graph, one has to
consider the result of breaking a framed $4$-graph at some edge
midpoint; this will result in two vertices of valency $1$. All
vertices of valency $4$ inherit the framing (the half-edges incident
to a four-valent vertex are split into two pairs of formally
opposite edges); one may also say that the newborn vertices of
valency $1$ are mapped to $\infty$ and $+\infty$ and the
intersection of the image of the graph with the exterior of some
large circle coincides with $Ox$, and the graph is oriented from
$-\infty$ to $+\infty$ along its unicursal component.
\end{dfn}

\begin{dfn}
 A {\em
long virtual knot} is an equivalence class of long virtual diagrams
by the same Reidemeister moves. It is allowed to apply these
Reidemeister moves only inside some prefixed large circle.
\end{dfn}

\subsection{Virtual Knots and Gauss Diagrams}

\begin{dfn}
With a virtual knot $K$ diagram one naturally associates (see
\cite{GPV}) a {\em Gauss diagram} $\Gamma(K)$, with all chords
endowed with arrows and signs; chords of the Gauss diagrams are in
one-one corresondence with {\em classical} crossings of the diagram;
the arrow is pointed from the preimage of the overcrossing arc to
the preimage of the undercrossing arc; the sign of the crossing
locally looking like $\skcrro$ is positive, and the sign of a
crossing locally looking like $\skcrlo$ is negative. With a (long)
virtual knot diagram $K$ one associates a (long) free knot whose
chord diagram is obtained from $\Gamma(K)$ by forgetting arrows and
signs. Note that the detour move does not change the Gauss diagram
corresponding to a virtual diagram. By construction, Reidemeister
moves on virtual knots generate the equivalence relations on free
knots described above
also called {\em Reidemeister moves}.
\end{dfn}

The Gauss diagram of classical and virtual trefoils are given below.

 \begin{figure}
 \centering\includegraphics[width=200pt]{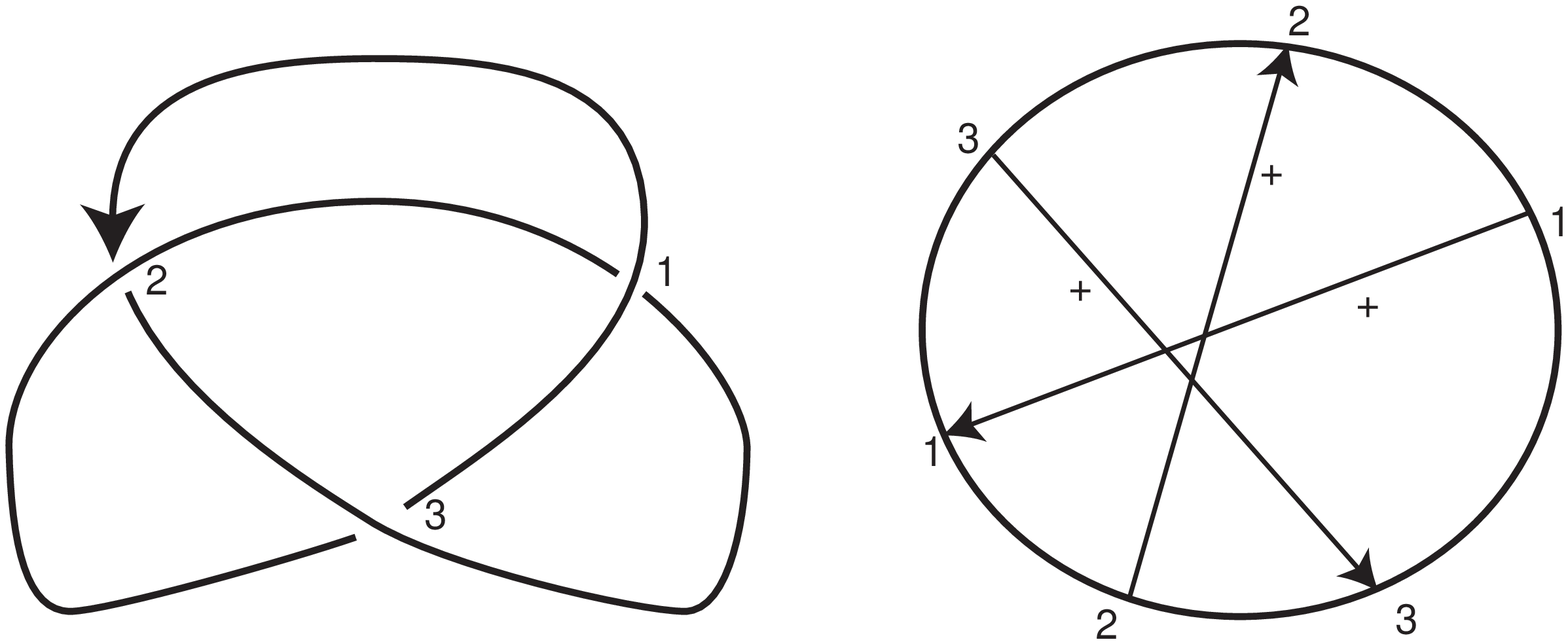}
  \caption{The Right Trefoil and Its Gauss Diagram}
  \label{Gaussdiagram}
 \end{figure}

 \begin{figure}
  \centering\includegraphics[width=200pt]{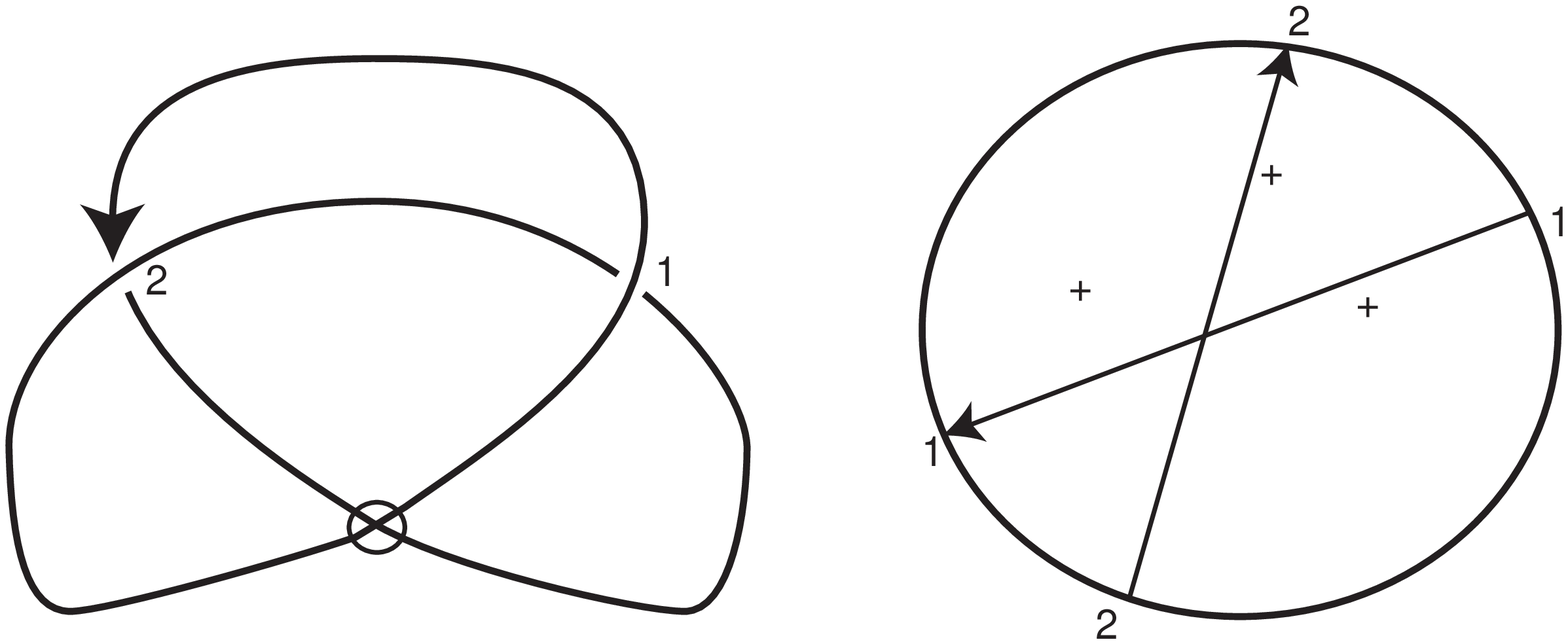}
\caption{The Virtual Trefoil and Its Gauss Diagram}
\label{VGaussdiagram}
 \end{figure}

\begin{figure}
\centering\includegraphics[width=120pt]{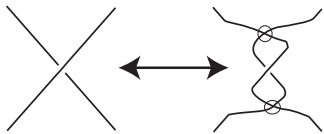} \caption{The
virtualization move} \label{virtua}
\end{figure}

By {\em virtualization move} for a classical crossing we mean the
local transformation of a (virtual) diagram, which inverts the arrow
direction for the arrow corresponding to this crossing, and
preserves the sign of the crossings (in the level of Gauss
diagrams). This map is well defined up to detour moves. In the
language of virtual diagrams the virtualization move looks as shown
in Fig. \ref{virtua}.

In \cite{Chr}, Chrisman proved that Goussarov-Poylak Viro
combinatorial formulae for virtual knots do not yield invariants
preserved by virtualization. The invariants of virtual knots
constructed in the present paper are all invariant under
virtualization.

\section{Vassiliev (finite-type) invariants. The main result}

Let $h$ be an invariant of (long) virtual knots. We say that  $h$ is
a {\em finite-type (Vassiliev)} invariant of order at most  $k$, if
for every virtual knot $K$ with $k+1$ fixed classical crossings the
following takes place $\sum_{s} (-1)^{\#(s)} h(K_{s})=0$, where the
sum is taken over all possible  $2^{k+1}$ knots $K_{s}$ coiniciding
with  $K$ everywhere except neighborhoods of given crossings, and at
the fixed $k+1$ places, the crossings $\skcrro$ or $\skcrlo$ are
chosen arbitrarily. Here the symbol $\#(s)$ stays for the number of
crossings of type  $\skcrlo$ among the chosen ones.

The other way of defining finite-type invariants uses {\em singular
knots} and {\em rigid vertices}, for more details see, \cite{KaV}. A
rigid vertex $\skcross$ is a formal linear combination
$\skcross=\skcrro-\skcrlo$. A {\em singular knot} of order $k$ is a
formal linear combination of $2^{k}$  knots which appear as
resolutions of $k$ formal singular crossings (rigid vertices). Every
knot invariant naturally extends to linear combinations of knots,
hence, to singular knots. An invariant is of finite type at most $k$
whenever its extension to singular knots of order $k+1$ vanishes.

Our next goal is to use a wider group than $G_{m}$ to construct
invariants of {\em virtual knots} (with some over/undercrossing
information). This goal can be achieved by extending the group
$G_{m}$.

Consider the group ${\tilde G}_{m}=\langle a'_{0},a''_{0},\dots,
a'_{m-1},a''_{m-1},a_{m}|a_{m}^{2}=e,(a'_{i})^{2}=(a''_{i})^{2},
i=0,\dots, m\!-\!1, a'_{i}a'_{j}=a'_{j}a''_{i},
a'_{i}a''_{j}=a''_{j}a''_{i},
a'_{j}a'_{i}=a''_{i}a'_{j},a''_{j}a'_{i}=a''_{i}a''_{j}, i<j;
a'_{i}a_{m}=a_{m}a''_{i}, i=0,\dots, m-1 \rangle$.

Let us fix a natural number $m$, and write $G=G_{m}$. It is clear
that the group $G$ defined above is obtained from ${\tilde G}_{m}$
by adding the following defining relations:
$(a'_{i})^{2}=1,i=0,\dots, m-1$. Let ${\cal G}=\Q{\tilde G}_{m}$ be
the group algebra of${\cal G}$, and let ${\cal G}_{k}$ denote the
quotient algebra of ${\cal G}$ by the following relations
$\prod_{j=1}^{k+1} ((a'_{n_{j}})^{2}-1)=0$, where $n_{j}$ stays for
any arbitrary set of numbers from $0$ to $m-1$. It is clear that
${\cal G}_{0}=\Q G$.

From the relations for the group $G$, it easily follows that every
element of ${\cal G}_{k}$ looking like  $\alpha_{1}\cdot A_{1}\cdot
\alpha_{2} A_{2}\cdots \alpha_{k+1} A_{k+1}\cdot \alpha_{k+2}$ for
arbitrary $\alpha$ and $A_{j}=(u-u^{-1})$, where $u$ stays for a
generator of the group $G$,  is  {\em equal to zero}. Indeed, it
suffices to note that the commutation relations for $u$ and $u^{-1}$
are similar, e.g. for $u=a'_{1}$ we have $a'_{i}u=a''_{1}$ and
$a'_{i}u^{-1}=(a''_{1})^{-1}$ for $i>1$. So, $(u-u^{-1})$ in this
case transforms into $a''_{1}-(a''_{1})^{-1}$. So, we can collect
all expressions of type $(u-u^{-1})$ together, and get $0$.

Let  $K$ be a long virtual knot, (see \cite{MaLong}). With it, we
associate an element  $\delta(K)\in {\tilde G}\subset {\cal G}$ as
follows. Take a Gauss diagram $\Gamma(K)$ of $K$, and  start writing
a word in generators of the group ${\tilde G}$ and its inverses
exactly in the way we were writing down the word $\gamma(D)$, with
the only difference that instead of each generator $a'_{j}$ or
$a''_{j}, j=0,\dots, m-1$ we shall write down either this generator
or its inverse depending on whether the crossing in question is
positive $(\skcrro)$ or negative $(\skcrlo)$. Denote the obtained
word by $\delta(K)$. Note that the word $\delta(K)$ by definition
does not depend on the direction of arrows in the Gauss diagram of
$K$.

The main result of the present paper is the following
\begin{thm}
The element $\delta(K)$ of the group $\tilde G$ is an invariant of
long virtual knots, which does not change under virtualization
(change of arrow direction on the Gauss diagram).

The conjugacy class of the element $\delta(K)$ in $\tilde G$ is an
invariant of (compact) virtual knot.

For every  $k$, the map $K\to \delta(K)\in {\cal G}_{k}$ is a
Vassiliev invariant of long virtual knots of order less than or
equal to  $k$.
\end{thm}

\begin{proof}
The first statement of the theorem easily follows from a comparison
of Reidemeister moves for virtual knots and the relations in the
group ${\tilde G}$; by construction, the invariant $\delta$ does not
depend on the arrow direction on the Gauss diagram: the latter
statement means the invariance under virutalization (see
\cite{FKM}). When passing from long virtual knots to compact virtual
knots, we allow the marked point to go along the Gauss diagram (or,
equivalently, we allow the infinity change); when moving the marked
point through a crossing we conjugate the correspondent element of
${\tilde G}$, which yields the second statement of the theorem.

Now, consider the alternating difference of the values of  $\delta$
for all long virtual knots corresponding to  $k+1$ choices crossing
types $\skcrro$ and $\skcrlo$ (or, equivalently, the value of
$\delta$ on the singular knot of order $k+1$). By construction of
 $\delta$, each singular crossing will contribute a factor of
 one of the following types: $(a'_{i}-(a'_{i})^{-1})$ or $(a''_{i}-(a''_{i})^{-1})$ or
 $(a_{m}-(a_{m})^{-1})$; the latter factor is zero. Taking into
 account the factorization relation defining  ${\tilde G}_{k}$,
 and the fact that the number of factors in our product is equal to $k+1$,
we conclude that the desired alternating sum in ${\tilde G}_{k}$ is
equal to zero.
\end{proof}

The invariants described in the present paper are constructive. The
group ${\tilde G}_{m}$ admits a simple description  similar to that
of the group $G$. The Cayley graph of the group ${\tilde G}_{m}$ is
the integer grid in  $\R^{2n+1}$; all coordinates of the vertices of
this graph are arbitrary integer numbers, except the last one, which
is equal to either zero or one. The right multiplication by $a_{m}$
changes the last coordinate, whence the right multiplication by
$a'_{i}$ or $a''_{i}$ corresponds to a shift in positive direction
along one of coordinates $x_{2i+1}$ or $x_{2i+2}$; here the choice
which coordinate changes depends on the parity of sum of all
coordinates from $2i+1$ to $m$.

\section{Values of the Invariants and Further Discussion}

\begin{figure}
\centering\includegraphics[width=200pt]{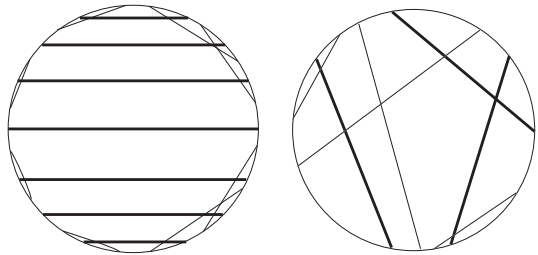} \caption{Two
non-slice free knots} \label{exs}
\end{figure}

First we note that even for free knots and even for the group
$G_{1}$, the values of the invariant $[\gamma]$ are very
interesting. In Fig. \ref{exs}, we give two knots where the two
coordinates of these values are equal to $16$ and $8$, respectively.

The knots above are not {\em slice} in the sense of paper
\cite{Ma4}. In fact, the invariant $[\gamma]$ provides a sliceness
obstruction for free knots, as shown in \cite{Ma4}.

It is not a difficult exercise to show that the $L(K)$ is divisible
by $4$ for every free knot. The divisibility of this invariant  is
conjectured by O.V.Manturov and will be proved elsewhere.

So, one can see that even the invariant of free knots and even for
the case of the group $G_{1}$ is highly non-trivial. Certainly, so
are invariants of long virtual knots.

In the present paper, we have constructed finite-type invariants of
{\em long} virtual knots. For compact virtual knots, this can be
done as well, but not in that elegant manner: one cannot further
deal with specific elements of ${\tilde G}_{m}$, but should rather
take conjugacy classes; the problem of extracting conjugacy classes
from one another and taking alternating sum looks cumbersome. We
shall touch on this problem in a separate publication.


\begin{thebibliography}{100}

\bibitem[AM]{AM} J.E. Andersen and J. Mattes,
Configuration space integrals and universal Vassiliev
invariants over closed surfaces, arXiv:q-alg/9704019.

\bibitem[BL]{BL} J.Birman, X-S.Lin (1993), Knot Polynomials
and Vassiliev's Invariants, {\em Inventiones Mathematicae}, {\bf
111}, P. 225-270.

\bibitem[BN]{BN} D.Bar-Natan, On the Vassiliev Knot Invariants, (2005), {\em
Topology}, {\bf 34}, pp. 423-475.


\bibitem[Chr]{Chr} M.Chrisman, Twist Lattices and the Jones-Kauffman Polynomial for Long Virtual
Knots. {\em J. Knot Theory \& Ramif.}, to appear.

\bibitem[FKM]{FKM} R.A.Fenn, L.H.Kauffman, V.O.Manturov (2005),
Virtual Knots: Unsolved Problems, {\em Fund. Math.},{\bf 188}, pp.
293-323.

\bibitem[Gib]{Gib} A.Gibson, Homotopy invariants of Gauss words, ArXiv:
Math.GT/0902.0062



\bibitem[GPV]{GPV} Goussarov M., Polyak M., and Viro O (2000),
Finite type invariants of classical and virtual knots, {\em
Topology}. {\bf 39}. P. 1045--1068.

\bibitem[Gou]{Gus} M.N. Goussarov (1991), Novaya forma polinoma Jones'a-Conway'a
dlya orientirovannyh zacepleniy (A new form of the Jones-Conway
polynomial for oriented links),  {\em Zap. Nauchn. Seminarov LOMI.}
{\bf 193}, Geometry and Topology, 1. P. 4-9.

\bibitem[HK]{HK} D.Hrencecin, L.H.Kauffman (2003), On filamentations and virtual knots, {\em Topology Appl.} {\bf
134}, pp. 23-52.

 \bibitem[Ka1]{KaV}
L.\,H.~Kauffman (1999), Virtual Knot Theory, {\em Eur. J.
Combinatorics}  {\bf 20} (7), P.\ 662--690.

\bibitem[Kon]{Kon}
M.\,Kontsevich (1993), Vassiliev's Knot Invariants, {\em Adv. Sov.
Math.}, {\bf 16} (2), P. 137-150.

\bibitem[Ma]{Ma} V.O.Manturov (2010), {\em Parity in Knot Theory},
Sbornik Mathematics., {\bf 201} (5), pp. 65-110.

\bibitem[Ma2]{MaLong} V.O.Manturov (2005), {\em On Long Virtual Knots},
Doklady Mathematics, {\bf 141} (5), С. 195-198.

\bibitem[Ma3]{Mybook} V.O.Manturov (2003), {\em Knot Theory},
Chapman and Hall/CRC., 416 pp.

\bibitem[Ma4]{Ma4} V.O.Manturov (2010), {\em Parity and Cobordisms of Free Knots},
arXiv.Math/GT:1001.2827

\bibitem[MM]{MM}
O.V.Manturov, V.O.Manturov, Free Knots and Groups (2010), {\em
Journal of Knot Theory and Its Ramifications}, {\bf 19}, (2)

\bibitem[MM2]{MM2} O.V.Manturov, V.O.Manturov, Svobodnye Uzly i
Gruppy (Free Knots and Groups), {\em Doklady Mathematics}, to
appear.

\bibitem[Tu]{Tu} V.G.Turaev,Topology of words, {\em Proc. Lond. Math. Soc.} (3) 95 (2007), no. 2, С. 360–412.

\bibitem[Vas]{Vas}
V.A.Vassiliev (1990), Cohomology of Knot Spaces, In: Theory of
Singularities and Its Applications, {\em Adv. Sov. Math.}, {\bf 1}
(23), P. 23-70

\end{thebibliography}
\end{document}